\documentclass[12pt]{article}
\usepackage{amsmath}
\usepackage{amsfonts}
\usepackage{amssymb}
\usepackage{amsthm}
\usepackage{geometry}                
\usepackage{graphicx}
\usepackage{epstopdf}
\begin{document}

\title{{\bf  Solutions of diophantine equations as periodic points of $p$-adic algebraic functions, I}}         
\author{Patrick Morton}        
\date{}          
\maketitle

\begin{abstract}
Solutions of the quartic Fermat equation in ring class fields of odd conductor over quadratic fields $K=\mathbb{Q}(\sqrt{-d})$ with $-d \equiv 1$ (mod $8$) are shown to be periodic points of a fixed algebraic function $T(z)$ defined on the punctured disk $0< |z|_2 \le \frac{1}{2}$ of the maximal unramified, algebraic extension $\textsf{K}_2$ of the $2$-adic field $\mathbb{Q}_2$.  All ring class fields of odd conductor over imaginary quadratic fields in which the prime $p=2$ splits are shown to be generated by complex periodic points of the algebraic function $T$, and conversely, all but two of the periodic points of $T$ generate ring class fields over suitable imaginary quadratic fields.  This gives a dynamical proof of a class number relation originally proved by Deuring.  It is conjectured that a similar situation holds for an arbitrary prime $p$ in place of $p=2$, where the case $p=3$ has been previously proved by the author, and the case $p=5$ will be handled in Part II.
\end{abstract}

\section{Introduction.}

In this paper and its sequel it will be shown that the periodic points of an algebraic function, suitably defined, have, in several particularly interesting cases, number theoretic significance. I shall primarily consider algebraic functions defined on subsets of $\mathfrak{p}$-adic fields.  \smallskip

An important problem in algebraic number theory is to classify the finite extensions $L$ of an algebraic number field $K$ for which $\textrm{Gal}(L/K)$ is abelian.  These are the {\it abelian extensions} of $K$, and for certain fields $K$ we have a good understanding of how to find explicit generators for these extensions.  For example, a famous theorem known as the Kronecker-Weber Theorem says that all abelian extensions of the rational field $K=\mathbb{Q}$ are subfields of cyclotomic fields $\mathbb{Q}(\zeta_f)$, where $\zeta_f$ is a primitive $f$-th root of unity with $f \ge 3$.  In the case that $K=\mathbb{Q}(\sqrt{-d})$ is an imaginary quadratic extension of $\mathbb{Q}$, the abelian extensions of $K$ are known to be subfields of {\it ray class fields}, where the latter fields are generated over $K$ by the coordinates of torsion points on elliptic curves $E$ whose endomorphism rings are isomorphic to certain subrings (orders) of $K$ (see [h1] and [si]).  Such a curve is said to have complex multiplication by the subring $\textsf{R} \subset R_K$, where $R_K$ is the ring of algebraic integers contained in $K$, if $\textrm{End}_{\overline{\mathbb{Q}}}(E) \cong \textsf{R}$ and $\mathbb{Z} \subsetneq \textsf{R}$.  There is an important subclass of these abelian extensions known as {\it ring class fields}, which are generated over $K$ by the $j$-invariants $j(E)$ of elliptic curves $E$ with complex multiplication by orders contained in $K$.  The properties of ring class fields are developed in the classical theory of complex multiplication, which is the main focus of the book by Cox [co1]. \medskip

In class field theory (see [ch], [d3], or [h1]), the ring class fields over $K$ are characterized as follows.  If $f$ is a positive integer, the ring class field (mod $f$) of $K=\mathbb{Q}(\sqrt{-d})$, denoted by $\Omega_f$, is the unique abelian extension of $K$ having the property that the prime ideals $\mathfrak{p}$ (not dividing $f$) of the ring of integers $R_K$ of $K$, which split completely into prime ideals of degree $1$ in the ring of integers $R_{\Omega_f}$ of $\Omega_f$, are exactly those $\mathfrak{p}$ for which $\mathfrak{p}=(\xi)$ is principal in $R_K$ with $\xi \equiv r$ (mod $f$) and $r \in \mathbb{Z}$.  It follows from class field theory that $\textrm{Gal}(\Omega_f/K) \cong A_f/P_f$, where $A_f$ is the group of fractional ideals of $K$ which are relatively prime to $f$ and $P_f$ is the subgroup of $A_f$ consisting of principal ideals of the form $(\xi)$ for numbers $\xi \equiv r$ (mod $f$) and $r \in \mathbb{Z}$.   The set of all such integers $\xi$ of $R_K$ is a ring $\textsf{R}_{-d}$, which gives rise to the name {\it ring class field}.  If $d_K$ is the discriminant of $K$, the integer $-d=d_Kf^2$ is called the discriminant of the ring (order) $\textsf{R}_{-d}$.  In [co1, pp. 190-192] it is shown that the subfields of the fields $\Omega_f$ are exactly the abelian extensions $L$ of $K$ for which $\textrm{Gal}(L/\mathbb{Q})$ is a generalized dihedral group. \medskip

Let $\textsf{K}_p$ be the maximal unramified, algebraic extension of the $p$-adic field $\mathbb{Q}_p$.  Call an imaginary quadratic field $K$ $p$-admissible, for a given prime $p \in \mathbb{Z}$, if $\displaystyle \left(\frac{d_K}{p}\right)=+1$, where $d_K$ is the discriminant of $K$, so that $p$ splits into two prime ideals in the ring of integers $R_K$.  If $K$ is $p$-admissible, then its discriminant is a square in $\mathbb{Q}_p$, and $K$ can therefore be embedded in $\mathbb{Q}_p$.  Moreover, if $p \nmid f$, then $\Omega_f/K$ is unramified at $p$ and can also be embedded in $\textsf{K}_p$.  My goal in this paper is to prove a special case of the following conjecture, which was stated in [m3].\medskip

\noindent {\bf Conjecture 1.} {\it Let $p$ be a fixed prime number.  There is an algebraic function $T_p(z)$, defined on a certain subset $\textsf{D}_p \subseteq \textsf{K}_p$ of the maximal unramified, algebraic extension of $\mathbb{Q}_p$, such that $T_p(\textsf{D}_p) \subseteq \textsf{D}_p$, with the following properties:} \smallskip

\noindent (a) {\it Any ring class field $\Omega_f \subset \textsf{K}_p$ of a $p$-admissible field $K \subset \mathbb{Q}_p$, whose conductor $f$ is relatively prime to $p$, is generated over $K$ by a periodic point $\xi$ of $T_p(z)$ contained in $\textsf{D}_p$;} \smallskip

\noindent (b) {\it All but finitely many periodic points $\xi$ of $T_p(z)$ contained in $\overline{\mathbb{Q}}_p$ generate ring class fields $\Omega_f=K(\xi)$ over some $p$-admissible quadratic field $K$.} \bigskip

The situation referred to in this conjecture is analogous to the fact that the fields $\mathbb{Q}(\zeta_f)$, where $\zeta_f$ is a primitive $f$-th root of unity and $(f,p)=1$, are generated over $\mathbb{Q}$ by periodic points of the map $F(z)=z^p$.  In fact, $\zeta_f$ is a periodic point of $F(z)$ with period $n$, where $n$ is the order of the prime $p$ modulo $f$.  Furthermore, the fields $\mathbb{Q}(\zeta_{p^k f})$ are generated over $\mathbb{Q}$ by {\it pre-periodic} points of $F(z)$, since $\zeta_{p^k f}$ is a root of $F^{k+n}(z)-F^k(z)=0$, for the same value of $n$.  \smallskip

An algebraic function $T_3(z)$ satisfying Conjecture 1 for the prime $p=3$ was given in [m3], namely
\begin{equation*}
T_3(z)=\frac{z^2}{3}(z^3-27)^{1/3}+\frac{z}{3}(z^3-27)^{2/3}+\frac{z^3}{3}-6, \ \ \textrm{for} \ z \in \textsf{K}_3, \ |z|_3 \ge 1,
\end{equation*}
where $T_3(z)$ is defined using the binomial series.  The periodic points of the function $T_3(z)$ in its $3$-adic domain were shown to be solutions of the cubic Fermat equation in ring class fields $\Omega_f$ over $3$-admissible quadratic fields $K=\mathbb{Q}(\sqrt{-d})$, whose conductors $f$ are prime to $3$.  Furthermore, every such $\Omega_f$ is generated over $\mathbb{Q}$ by one of these periodic points.  \smallskip

In this paper I will show that a certain $2$-adic branch of the function
\begin{equation*}
T(z)=\frac{\sqrt[4]{1-z^4}+1}{\sqrt[4]{1-z^4}-1}=1-\frac{2}{z^4}\left(1+(1-z^4)^{1/4}+(1-z^4)^{1/2}+(1-z^4)^{3/4}\right)
\end{equation*}
satisfies the statement of the above conjecture for the prime $p=2$.  I will show that all of the periodic points of $T(z)$ in its $2$-adic domain $\textsf{D}_2=\{z: 0 < |z|_2 \le \frac{1}{2} \} \subset \textsf{K}_2$ are solutions of the quartic Fermat equation in ring class fields of $2$-admissible quadratic fields.  These solutions have been given in [lm] as follows.  Though the precise formulas are not necessary for the proofs in this paper, it is worth noting that these solutions can be represented in terms of modular functions.  \smallskip

Let $\eta(\tau)$ be the Dedekind $\eta$-function [co1, p. 256].  The Schl\"afli functions $\mathfrak{f}(\tau), \mathfrak{f}_1(\tau)$, $\mathfrak{f}_2(\tau)$ (see [sch, p. 148] or [co1, p. 256]) are defined to be:

\begin{equation*}
\mathfrak{f}(\tau)=e^{-\frac{\pi i}{24}} \frac{\eta\left( \frac{\tau + 1}{2} \right) }{\eta(\tau)}, \quad \mathfrak{f}_1(\tau)= \frac{\eta\left( \frac{\tau }{2} \right) }{\eta(\tau)},\quad \mathfrak{f}_2(\tau)= \sqrt{2} \hspace{.05 in} \frac{\eta(2\tau) }{\eta(\tau)}.
\end{equation*}
These functions have the infinite product representations

\begin{equation*}
\mathfrak{f}(\tau)=q^{-\frac{1}{48}} \prod_{n=1}^\infty{(1+q^{n-\frac{1}{2}})},\quad \mathfrak{f}_1(\tau)=q^{-\frac{1}{48}} \prod_{n=1}^\infty{(1-q^{n-\frac{1}{2}})},
\end{equation*}

\begin{equation*}
\mathfrak{f}_2(\tau)=\sqrt{2} \hspace{.05 in} q^{\frac{1}{24}} \prod_{n=1}^\infty{(1+q^n)}, \quad q = e^{2 \pi i \tau},
\end{equation*}
convergent on the upper half-plane $\mathbb{H}$.  Let $K=\mathbb{Q}(\sqrt{-d})$ be a 2-admissible quadratic field, where $-d \equiv 1$ (mod $8$) is the discriminant of the order $\textsf{R}_{-d}$ in $K$, with conductor $f$, satisfying $-d=d_K f^2$.  Further, let $w \in K$ be defined by
\begin{equation*}
w=\frac{v+\sqrt{-d}}{2}, \quad v^2 \equiv -d \ (\textrm{mod} \ 16), \ \ v=1 \ \textrm{or} \ 3,
\end{equation*}
and set
\begin{equation*}
a \equiv
\begin{cases}
\frac{-3d+5}{\ 16} \ (\textrm{mod} \ 4), &\textrm{if} \ v=3 \ \textrm{and} \ d \equiv 7 \ (\textrm{mod} \ 16),\cr
	\frac{-d+31}{16} \ (\textrm{mod} \ 4), &\textrm{if} \ v=1 \ \textrm{and} \ d \equiv 15 \ (\textrm{mod} \ 16).
\end{cases}
\end{equation*}
Then the numbers
\begin{equation}
\pi_d = i^a\frac{\mathfrak{f}_2(w/2)^2}{\mathfrak{f}(w/2)^2}, \ \ \xi_d=\frac{\beta}{2}=i^{-v}\frac{\mathfrak{f}_1(w/2)^2}{\mathfrak{f}(w/2)^2}
\end{equation}
lie in the ring class field $\Omega_f$ of conductor $f$ over $K$, and satisfy
\begin{equation*}
\pi_d^4+\xi_d^4=1.
\end{equation*}
(See [lm, Sec. 10].)  The numbers $\pi_d$ and $\xi_d$ are conjugate algebraic integers over $\mathbb{Q}$ and $\Omega_f$ is generated over $\mathbb{Q}$ by either of them.  Furthermore, if $\wp_2=(2,w)$ is one of the prime ideal divisors of $2$ in $K$, then with $(2)=2R_K=\wp_2 \wp_2'$, we have
\begin{equation*}
(\pi_d)=\pi_d R_{\Omega_f} = \wp_2 R_{\Omega_f}, \ \ (\xi_d)=\xi_d R_{\Omega_f} = \wp_2'R_{\Omega_f}, \ \textrm{in} \ \Omega_f,
\end{equation*}
where $R_L$ denotes the ring of algebraic integers in the field $L$.  In other words, $\pi_d$ and $\xi_d$ are principal ideal generators in $R_{\Omega_f}$ of the prime ideal divisors of $2$ in $R_K$, when those ideals are extended to the larger ring $R_{\Omega_f}$.  \smallskip

Denote by $b_d(x)$ the minimal polynomial over $\mathbb{Q}$ of the numbers $\pi_d$ and $\xi_d$.  Then $b_d(x)$ is a normal polynomial over $\mathbb{Q}$ (meaning that one of its roots generates a normal extension of $\mathbb{Q}$) and
\begin{equation*}
\textrm{deg}(b_d(x))=2h(-d),
\end{equation*} 
where $h(-d)=|A_f/P_f|$ is the class number of the order $\textsf{R}_{-d}$, i.e. the number of elements of the ideal class group of $\textsf{R}_{-d}$.  See [co1, pp. 132-148]; and see Section 6 for some examples of these polynomials. \smallskip

To explicitly define the branch of $T(z)$ that we will be considering, let
\begin{equation*}
T_1(z)=\frac{2z^4-4-4\sqrt{1-z^4}}{z^4}, \ \ T_2(z)=\frac{z}{2}-\frac{z}{2}\sqrt{1-\frac{4}{z^2}},
\end{equation*}
where the square-roots are defined $2$-adically by the binomial series.  We have: \bigskip

\noindent {\bf Theorem 1.}  {\it a) The function $T(z)=T_2 \circ T_1(z)$ maps the set
\begin{equation*}
\textsf{D}_2=\{z: 0 < |z|_2 \le \frac{1}{2} \} \subset \textsf{K}_2
\end{equation*}
to itself.  \smallskip

b) The periodic points of $T(z)$ in $\textsf{D}_2$ are the roots $\xi_d$ of the polynomials $b_d(x)$, as $-d$ varies over quadratic discriminants $\equiv 1$ (mod $8$), along with the conjugates of $\xi_d$ over $K=\mathbb{Q}(\sqrt{-d})$, under the natural embedding of $\Omega_f$ in its completion $(\Omega_f)_\mathfrak{p} \subset \textsf{K}_2$, for a prime ideal $\mathfrak{p}$ of $R_{\Omega_f}$ which divides $\wp_2'$.  
\smallskip

c) The number of periodic points of $T(z)$ in the domain $\textsf{D}_2$ with minimal period $n$ is given by
\begin{equation*}
\sum_{-d \in \mathfrak{D}_n}{h(-d)}=nN_4(n)=\sum_{k \vert n}{\mu (n/k)2^{2k}}, \quad n > 1.
\end{equation*}
Here $\mathfrak{D}_n$ is the set of discriminants $-d \equiv 1$ (mod $8$) for which the square of the corresponding Frobenius auotmorphism $\tau=\left(\frac{\Omega_f/K}{\wp_2}\right)$ has order $n$ in $Gal(\Omega_f/K)$, and $\mu$ is the M\"obius $\mu$-function.  For $n=1$, the number of fixed points of $T(z)$ in $\textsf{D}_2$ is}
\begin{equation*}
\sum_{-d \in \mathfrak{D}_1}{h(-d)}=h(-7)+h(-15)=4-1.
\end{equation*}
\medskip

Thus, our analysis gives a dynamical interpretation of the class number formula occurring in part c), which is equivalent to a special case of a class number formula of Deuring [d1], [d2, Eq. (1)].  As in [m3], the proof implies a similar statement about the periodic points of the multi-valued function $T(z)$ on either of the fields $\overline{\mathbb{Q}}_2$ or $\mathbb{C}$.  For this it is necessary to define what we mean by a periodic point of a multi-valued algebraic function. \smallskip

Let $f(z)$ be any algebraic function defined over a given field $F$, so that $f(z)$ lies in the algebraic closure $\overline{F(z)}$ of $F(z)$, and let $g(z,w) \in F[z,w]$ be the minimal polynomial of $w=f(z)$ over $F(z)$.  Define a periodic point $a$ over $F$ of the algebraic function $f(z)$ to be any number $a \in F$ for which there exist $a_1, a_2, \cdots, a_{n-1} \in F$ for which
\begin{equation*}
g(a,a_1)=g(a_1,a_2)=\cdots =g(a_{n-2},a_{n-1})=g(a_{n-1},a)=0.
\end{equation*}
By cyclically permuting the equations in the definition it is clear that all the numbers $a_i$ are also periodic points of $f(z)$ of period $n$.  Thus, when writing $f(a_{i-1})=a_i$, each individual element $a_i=f_i(a_{i-1})$ will be defined using one particular branch $f_i(z)$ of $f(z)$, for $1 \le i \le n$ (taking $a_0=a_n=a$), and different branches $f_i, f_j$ may or may not coincide.  With this definition we have the following. \bigskip

\noindent {\bf Theorem 2.} {\it The set of periodic points of the multi-valued function $T(z)$ on any of the fields $\mathcal{K}=\textsf{K}_2$, $\overline{\mathbb{Q}}_2$ or $\mathbb{C}$ coincides with the set
$$\mathcal{S(K)}=\{0,-1\} \cup \{\xi \in \mathcal{K}: (\exists n \ge 1)(\exists \ (-d) \in \mathfrak{D}_n ) \ \textrm{s.t.} \ b_d(\xi)=0 \}.$$
Thus, all the periodic points of $T(z)$ distinct from $0$ and $-1$ in any of these fields generate ring class fields over $2$-admissible quadratic extensions of $\mathbb{Q}$, and give solutions of the quartic Fermat equation.}  \medskip

In part II of this paper, I shall verify the above conjecture for the prime $p=5$, by considering solutions of the diophantine equation
\begin{equation*}
\varepsilon^5X^5+\varepsilon^5Y^5=1-X^5Y^5, \quad \varepsilon=\frac{1+\sqrt{5}}{2},
\end{equation*}
in certain class fields of $5$-admissible quadratic fields. \smallskip

The following conjecture is also stated in [m3].  \bigskip

\noindent {\bf Conjecture 2.}  {\it Any ring class field of a $p$-admissible quadratic field $K=\mathbb{Q}(\sqrt{-d}) \subset \overline{\mathbb{Q}}_p$, whose conductor is divisible by $p$, is generated over $K$ by some pre-periodic point of the multi-valued function $T_p(z)$ contained in the algebraic closure $\overline{\mathbb{Q}}_p$}. \bigskip

This statement was proved for $p=3$ and the above function $T_3(z)$ in [m3] and will be proved for $p=2$ and the $2$-adic function $T(z)$ elsewhere. (See [cn].) \smallskip

\section{The quartic Fermat equation.}

The numbers $\pi_d$ and $\xi_d$ defined in (1) were shown in [lm] to be algebraic conjugates of each other over $\mathbb{Q}$.  This fact was deduced from the relationship
$$\pi_d^{\tau^2}=\frac{\xi_d+1}{\xi_d-1},$$
where $\tau$ is a certain automorphism in the Galois group of $\Omega_f/K$, uniquely defined by the condition that
$$\alpha^\tau \equiv \alpha^2 \ (\textrm{mod} \ \wp_2),$$
for all elements $\alpha$ of the ring of integers $\Omega_f$, $R_{\Omega_f}$. Actually, this congruence holds for all $\alpha \in \Omega_f$ whose denominators are relatively prime to $\wp_2$ -- these are the elements of $\Omega_f$ which are integral for $\wp_2$.  This automorphism is denoted by
$$\tau=\left(\frac{\Omega_f/K}{\wp_2}\right),$$
and is called the Frobenius automorphism for the prime ideal $\wp_2$ of $R_K$.  An automorphism of $\textrm{Gal}(\Omega_f/K)$ can be assigned to any prime ideal $\mathfrak{p}$ in $R_K$ which is relatively prime to $f$ (and therefore unramified in $\Omega_f$), satisfying
$$\alpha^\sigma \equiv \alpha^{\textrm{Norm}(\mathfrak{p})} \ (\textrm{mod} \ \mathfrak{p}), \ \ \alpha \in R_{\Omega_f}, \ \ \sigma=\left(\frac{\Omega_f/K}{\mathfrak{p}}\right),$$
where $\textrm{Norm}(\mathfrak{p})=|R_K/\mathfrak{p}|$ is the absolute norm of $\mathfrak{p}$.  (See [ch], [co1], or [h2].)  Recall that $f$ is the positive integer for which $K=\mathbb{Q}(\sqrt{-d})$ and $-d=d_Kf^2$, where $d_K$ is the discriminant of $K/\mathbb{Q}$.  Although the square-roots of the numbers $-d_Kf^2$ all generate the same quadratic field $K$, the degrees of the numbers $\pi_d$ and  $\xi_d$ and the field they generate over $\mathbb{Q}$ depend strongly on the parameter $f$.  We always assume $-d \equiv 1$ (mod $8$), so that $d_K$ and $f$ are odd integers.  \medskip

Replacing $x$ by $(x+1)/(x-1)$ in the Fermat equation $x^4+y^4=1$ leads to the curve $f(x,y)=0$ defined by the equation
\begin{equation}
f(x,y)=y^4(x-1)^4+8x(x^2+1).
\end{equation}
Writing $\pi=\pi_d, \xi=\xi_d$, the relation $(\pi^{\tau^2})^4+(\xi^{\tau^2})^4=1$ yields
\begin{equation}
f(\xi,\xi^{\tau^2})=0, \ \ \xi=\frac{\beta}{2}.
\end{equation}
It follows that $\xi^{\tau^2}$ can be considered as one of the values of the algebraic function

\begin{equation*}
y=S(x)=\sqrt[4]{\frac{-8x(x^2+1)}{(x-1)^4}}=\sqrt[4]{1-\left(\frac{x+1}{x-1}\right)^4}
\end{equation*}
at $x=\xi$.  It is natural to try to expand $S(x)$ as follows:

\begin{equation*}
S(x)=1+\sum_{k=1}^\infty{(-1)^k {\frac{1}{4} \atopwithdelims ( ) k} \left(\frac{1+x}{1-x} \right)^{4k}}.
\end{equation*}
Unfortunately, this cannot be expressed as a convergent $2$-adic series in powers of $x$, since
$$S(0)=1+\sum_{k=1}^\infty{(-1)^k {\frac{1}{4} \atopwithdelims ( ) k}}$$
does not even converge ($2$-adically).  Instead, we apply $\tau^{-2}$ to $f(\xi,\xi^{\tau^2})=0$, obtaining $f(\xi^{\tau^{-2}},\xi)=0$, and we consider $\xi^{\tau^{-2}}$ as one of the values of the inverse algebraic function
\begin{equation}
x=T(y)=\frac{\sqrt[4]{1-y^4}+1}{\sqrt[4]{1-y^4}-1}, \ \ f(x,y)=0,
\end{equation}
evaluated at $y=\xi$.  \medskip

We first find an expression for a particular 2-adic branch of the function $T(y)$.  Expanding and dividing $f(x,y)$ by $y^4$ gives
\begin{align*}
\frac{f(x,y)}{y^4}= & \ x^4+\frac{8-4y^4}{y^4}x^3+6x^2+\frac{8-4y^4}{y^4}x+1\\
= & \ x^4+tx^3+6x^2+tx+1, \ \ t=\frac{8-4y^4}{y^4}.
\end{align*}
Hence,
\begin{align*}
\frac{f(x,y)}{x^2y^4}= & \ \left(x^2+\frac{1}{x^2}\right)+t\left(x+\frac{1}{x}\right)+6\\
= & \ z^2+tz+4, \ \ z=x+\frac{1}{x}.
\end{align*}
Thus we have
$$z=\frac{-t \pm \sqrt{t^2-16}}{2}=\frac{2y^4-4 \pm 4\sqrt{1-y^4}}{y^4}.$$
We define
\begin{equation}
T_1(y)=\frac{2y^4-4 - 4\sqrt{1-y^4}}{y^4}.
\end{equation}
This function can be expanded into a $2$-adic Laurent series in $y$, as follows.
\begin{align*}
T_1(y)&=2-\frac{4}{y^4}-\frac{4}{y^4} \sum_{n=0}^\infty{{\frac{1}{2} \atopwithdelims ( ) n} (-1)^n y^{4n}}\\
&=2-\frac{4}{y^4}-\frac{4}{y^4} (1-\frac{1}{2}y^4-\frac{1}{8} y^8 - \cdots)\\
&=\frac{-8}{y^4}+4+4 \sum_{n=2}^\infty{{\frac{1}{2} \atopwithdelims ( ) n} (-1)^{n+1} y^{4n-4}}.
\end{align*}
It is not hard to verify that the series for $T_1(y)$ converges $2$-adically for $0<|y|_2 \le \frac{1}{2}$.  To see this, set $y=2y_1$.  With this substitution, the series becomes
\begin{equation}
T_1(2y_1)+\frac{1}{2y_1^4}-4 = \sum_{n=2}^\infty{2^{4n-2}{\frac{1}{2} \atopwithdelims ( ) n} (-1)^{n+1} y_1^{4n-4}} =  \sum_{n=2}^\infty{2^{2n-1}C_{n-1} y_1^{4n-4}},
\end{equation}
where $C_{n-1}=(-1)^{n+1}2^{2n-1}{\frac{1}{2} \atopwithdelims ( ) n} \in \mathbb{Z}$ is the Catalan number.  Hence, the coefficient of $y_1$ in the series (6) is divisible by $2^{2n-1}$, and the series is therefore convergent for $|y_1| \le 1$.  This proves the above claim.  Moreover, the infinite series in (6) represents a $2$-adic integer for $|y|_2 \le \frac{1}{2}$, so it is clear that
\begin{equation}
|T_1(y)|_2 \ge 2, \ \ \textrm{if} \ 0 < |y|_2 \le \frac{1}{2},
\end{equation}
because of the leading term $\frac{-8}{y^4}$.  The second solution of $z^2+tz+4=0$ is then
$$-t-T_1(y)=-4 \sum_{n=2}^\infty{{\frac{1}{2} \atopwithdelims ( ) n} (-1)^{n+1} y^{4n-4}}=\frac{4}{T_1(y)},$$
which is a 2-adic integer.  \medskip

Solving the equation $x^2-zx+1=0$ for $x$ gives
$$x=\frac{z \pm \sqrt{z^2-4}}{2}.$$
Now we set
\begin{align*}
T_2(z)&=\frac{z}{2}-\frac{z}{2} \sqrt{1-\frac{4}{z^2}}=\frac{z}{2}-\frac{z}{2} \sum_{n=0}^\infty{(-1)^n {\frac{1}{2} \atopwithdelims ( ) n}\frac{2^{2n}}{z^{2n}}}\\
&=\sum_{n=1}^\infty{(-1)^{n+1} {\frac{1}{2} \atopwithdelims ( ) n}\frac{2^{2n-1}}{z^{2n-1}}}=\sum_{n=1}^\infty{\frac{C_{n-1}}{z^{2n-1}}}\\
&=\frac{1}{z}+\frac{1}{z^3}+\frac{2}{z^5}+\frac{5}{z^7}+\frac{14}{z^9}+ \frac{42}{z^{11}}+\cdots,
\end{align*}
which is convergent for $|z|_2 \ge 2$, as above.  It is clear from this series expansion that
\begin{equation}
0<|T_2(z)|_2 \le \frac{1}{2} \ \ \textrm{for} \ |z|_2 \ge 2,
\end{equation}
since $\frac{4}{z^2}\neq 0$.  The second solution of $x^2-zx+1=0$ is then $z-T_2(z)=\frac{1}{T_2(z)}$, which is {\it not} a 2-adic integer.  \medskip

By the above arguments, setting $z=T_1(y)$ gives the solution
\begin{equation}
x=T_2(z)=T_2(T_1(y)), \ \ \textrm{for} \ 0 < |y|_2 \le \frac{1}{2},
\end{equation}
of $f(x,y)=0$.  By (7) and (8), the function
$$T=T_2 \circ T_1$$
maps the region $0<|y|_2 \le \frac{1}{2}$ of $\textsf{K}_2$ into itself.  It is clear that this is also true of the region $|y|_2=\frac{1}{2}$. This is the branch of $T$ which we will use throughout our discussion.  To summarize, we have: \medskip

\noindent {\bf Proposition 3.} {\it The algebraic function $T(y)=T_2(T_1(y))$, where
\begin{align*}
T_1(y) &= \frac{-8}{y^4}+4+4 \sum_{n=2}^\infty{{\frac{1}{2} \atopwithdelims ( ) n} (-1)^{n+1} y^{4n-4}},\\
T_2(z) &= \sum_{n=1}^\infty{(-1)^{n+1} {\frac{1}{2} \atopwithdelims ( ) n}\frac{2^{2n-1}}{z^{2n-1}}},
\end{align*}
is defined on the punctured disk
$$\textsf{D}_2=\{y \in \textsf{K}_2: \ 0 < |y|_2 \le \frac{1}{2} \}$$
in the field $\textsf{K}_2$, and maps $\textsf{D}_2$ to itself.  For any $y \in \textsf{D}_2$, we have $f(T(y),y)=0$.} \bigskip

We now prove the following theorem. \bigskip

\noindent {\bf Theorem 4.} {\it Let $(\pi, \xi)$ be any solution of $X^4+Y^4=1$ in the ring class field $\Omega_f$ of odd conductor $f$ over $K=\mathbb{Q}(\sqrt{-d})$ which is conjugate over $K$ to the solution (1). Then under the embedding of $\Omega_f$ in the maximal unramified extension $\textsf{K}_2$ of the 2-adic field $\mathbb{Q}_2$ given by $\Omega_f \rightarrow \left(\Omega_f\right)_\mathfrak{p}$, where $\mathfrak{p}$ is a prime divisor of $\wp_2'$ in $R_{\Omega_f}$, we have
$$\xi^{\tau^{-2}}=T(\xi), \ \textrm{with} \ \tau^{-1}=\left(\frac{\Omega_f/K}{\wp_2'}\right),$$
where $T(y)$ is the $2$-adic algebraic function from Proposition 3.  Thus, $\xi \rightarrow T(\xi)$ is a lift of the square of the Frobenius automorphism corresponding to $\wp_2'$ on $\Omega_f/K$.} \bigskip

\noindent {\it Proof.} The Galois group $\textrm{Gal}(\Omega_f/K)$ is a generalized dihedral group (see [co1, pp. 190-191]), so the automorphism $\displaystyle \tau=\left(\frac{\Omega_f/K}{\wp_2}\right)$ (applied exponentially) satisfies
$$ \tau^{-1}=\phi^{-1} \tau \phi = \left(\frac{\Omega_f/K}{\wp_2'}\right),$$ 
(see [co1, p. 107]) where $\phi$ is an automorphism of $\Omega_f$ which restricts to the nontrivial automorphism of $K$, sending $\wp_2$ to its conjugate ideal $\wp_2'$.  Hence, we know that
$$\left(\frac{\xi}{2} \right)^{\tau^{-2}} \equiv \left(\frac{\xi}{2} \right)^4 \ (\textrm{mod} \ \wp_2') \ \ \textrm{in} \ \Omega_f.$$
Embedding $\Omega_f$ into $\textsf{K}_2$ by completing at a prime $\mathfrak{p}$ of $\Omega_f$ lying over $\wp_2'$, we obtain that the images of $\xi, \xi^{\tau^{-2}}$, which we denote by the same symbols, satisfy
\begin{equation*}
\left(\frac{\xi}{2} \right)^{\tau^{-2}} \equiv \left(\frac{\xi}{2} \right)^4 \ (\textrm{mod} \ 2) \ \ \textrm{in} \ \left(\Omega_f\right)_\mathfrak{p} \subset \textsf{K}_2,
\end{equation*}
and, since both sides of this congruence are units for $\wp_2'$, that
$$\frac{2^3 \xi^{\tau^{-2}}}{\xi^4} \equiv 1 \ (\textrm{mod} \ 2) \ \ \textrm{in} \ \left(\Omega_f\right)_\mathfrak{p} \subset \textsf{K}_2.$$
Now we have from (7) and the series for $T_2(z)$ that $T_2(T_1(y)) \equiv \frac{1}{T_1(y)}$ (mod $2^3$), so
\begin{equation}
\frac{2^3T(\xi)}{\xi^4}=\frac{2^3 T_2(T_1(\xi))}{\xi^4} \equiv \frac{2^3}{\xi^4T_1(\xi)} \equiv -1 \equiv 1 \ (\textrm{mod} \ 2), \ \textrm{for} \ 0 < |\xi|_2 \le \frac{1}{2}.
\end{equation}
It follows that
$$\frac{\xi^{\tau^{-2}}}{T(\xi)} =\eta^{-1} \equiv 1 \ (\textrm{mod} \ 2),$$
and therefore $T(\xi)=\eta \xi^{\tau^{-2}}$, where $\eta$ is a 2-adic unit.  But $T(\xi)$ and $\xi^{\tau^{-2}}$ are both roots of $f(x,\xi)=0$ in $\textsf{K}_2$.  From the above argument we know there is a second root of $f(x,\xi)=0$ in $\textsf{K}_2$ given by $T_1(\xi)-T_2(T_1(\xi))=T_1(\xi)-T(\xi)$, which is not a 2-adic integer, by (7), since $T(\xi) \in \textsf{D}_2$ by Proposition 3.  (Recall that $(\xi)=\wp_2'$ in $R_{\Omega_f}$, so that $|\xi |_2=\frac{1}{2}$ in $\textsf{K}_2$.)  Thus, $T(\xi)$ is distinct from this root.
\smallskip

Now I claim that the polynomial
$$g(x)=\frac{f(x,\xi)}{\xi^4}=x^4+tx^3+6x^2+tx+1, \ \ t=\frac{8-4\xi^4}{\xi^4},$$
has at most two roots in $\textsf{K}_2$.  To see this, note that the Ferrari cubic resolvent of $g(x)$ [co2, pp. 358-359], whose roots are rational expressions over $\mathbb{Q}_2(\xi)$ in the roots of $g(x)$, is
$$r(y)=y^3-6y^2+(t^2-4)y-2t^2+24=(y-2)(y^2-4y+t^2-12),$$
where the discriminant of the quadratic factor is given by
$$\delta=-4(t^2-16)=\frac{256(\xi^4-1)}{\xi^8}.$$
We have $1-\xi^4\equiv 1$ (mod $16$) since $|\xi |_2 = \frac{1}{2}$, so Hensel's Lemma implies that $\delta=-\mu^2$ for some $\mu \in \textsf{K}_2$. Therefore, $\sqrt{\delta} \notin \textsf{K}_2$, since $\mathbb{Q}_2(\sqrt{-1})$ is a ramified extension, and the resolvent $r(y)$ has exactly one root in $\textsf{K}_2$.  This shows that the polynomial $g(x)$ has exactly two roots in $\textsf{K}_2$ and that  $T(\xi)=\xi^{\tau^{-2}}$.  \smallskip

It is clear that the above discussion also holds for any conjugate of $\xi=\xi_d$ over $K=\mathbb{Q}(\sqrt{-d})$, since the ideal $\wp_2'$ is fixed by the elements of $\textrm{Gal}(\Omega_f/K)$, and since this Galois group is abelian.    $\square$  \bigskip

We use Theorem 4 to prove \bigskip

\noindent {\bf Theorem 5.}  {\it With notation as in Theorem 4, $\xi$ is a periodic point of the algebraic function $T(y)$ on the domain $\textsf{D}_2:=\{y : |y|_2 \le \frac{1}{2}\} \subset \textsf{K}_2$, whose period $n$ is equal to the order of the automorphism $\tau^{-2}$ in $Gal(\Omega_f/K)$.}  \medskip

\noindent {\it Proof.}  This follows from the fact that $\tau^{-2}$, as an automorphism on the completion $\left(\Omega_f\right)_\mathfrak{p}$ fixing the prime ideal $\wp_2'\mathbb{Z}_2=2\mathbb{Z}_2$ of $(R_K)_{\wp_2'}=\mathbb{Z}_2$, satisfies
\begin{equation*}
T(z)^{\tau^{-2}}=T(z^{\tau^{-2}}), \ \ \textrm{for} \ z \in \left(\Omega_f\right)_\mathfrak{p} \cap \textsf{D}_2,
\end{equation*}
since the coefficients of $T_1(2y_1)-\frac{1}{2y_1^4}$ (see (6)) and $T_2(z)$ lie in $\mathbb{Z}$.  Therefore,
$$T^2(\xi)=T(T(\xi))=T(\xi^{\tau^{-2}})=T(\xi)^{\tau^{-2}}=\xi^{\tau^{-4}},$$
and more generally, $T^k(\xi)=\xi^{\tau^{-2k}}, \ k \ge 1$.  Since $\xi$ generates $\Omega_f$ over $K$, we have $\xi^{\tau^{-2k}} \neq \xi$ for $k < n$.  Hence, $T^n(\xi)=\xi^{\tau^{-2n}}=\xi$, which shows that $\xi$ is a periodic point of $T$ with minimal period $n$.   $\square$ \bigskip
 
This proves part (a) of Conjecture 1 of the Introduction, since every ring class field $\Omega_f$ of odd conductor over the $2$-admissible field $K$ is generated by the coordinates of a solution of the quartic Fermat equation. \smallskip

We would now like to prove the converse; namely, that any periodic point of $T$ on the domain $\textsf{D}_2$ comes from one of the solutions $(\pi, \xi)$ in some ring class field $\Omega_f$ over $K=\mathbb{Q}(\sqrt{-d})$, with $-d \equiv 1$ (mod $8$).  \medskip

\section{Iterated resultants.}

Define the following iterated resultants, as in [m3].  Set $R^{(1)}(x,x_1)=f(x,x_1)$,
$$R^{(2)}(x,x_2)=Res_{x_1}(f(x,x_1),f(x_1,x_2)),$$
and recursively define
\begin{equation}
R^{(k)}(x,x_k)=Res_{x_{k-1}}(R^{(k-1)}(x,x_{k-1}),f(x_{k-1},x_k)), \quad k \ge 3.
\end{equation}
Then we set $x_n=x$ in $R^{(n)}(x,x_n)$ to obtain $R_n(x)$:
$$R_n(x) = R^{(n)}(x,x), \ \ n \ge 1.$$ \smallskip
\noindent From this definition it is easy to see that the roots of $R_n(x)$ are exactly the $a$'s for which there exist common solutions of the equations
\begin{equation}
f(a,a_1)=0, \quad f(a_1,a_2)=0, \quad \cdots \quad f(a_{n-1},a)=0.
\end{equation}
In particular, (12) holds for
$$a=\xi=T^n(\xi), \ a_{n-1}=T(\xi), \ a_{n-2}=T(a_{n-1})=T^2(\xi), \cdots, a_1=T(a_2)=T^{n-1}(\xi),$$
by Proposition 3 and Theorem 5, so that $T^k(\xi)$ is a root of $R_n(x)$ for any $k$ with $0 \le k \le n$.  It is straightforward to show by induction that
$$R^{(n)}(x,x_n) \equiv x_n^{4^n}(x+1)^{4^n} \ \ (\textrm{mod} \ 2),$$
and therefore
$$R_n(x) \equiv x^{4^n}(x+1)^{4^n} \ \ (\textrm{mod} \ 2).$$
In the following lemma, we show that $R_n(x)$ is monic and has degree $2 \cdot 4^n$. \bigskip

\noindent {\bf Lemma.} a) {\it For $n \ge 2$, $R^{(n)}(x,x_n)=A_n(x)x_n^{4^n}+S_n(x,x_n)$, where $A_n(x) \in \mathbb{Z}[x]$ is a monic polynomial satisfying $\textrm{deg}(A_n(x))=4^n$, and}
$$\textrm{deg}_{x_n}(S_n(x,x_n)) \le 4^n-4, \ \ \textrm{deg}_{x}(S_n(x,x_n)) \le 4^{n}-1.$$

b) {\it $\textrm{deg}(R_n(x))=2 \cdot 4^n$, and the leading coefficient of $R_n(x)$ is $1$}. \medskip

\noindent {\it Proof.} a) The assertion is obvious for $n=1$ by (2).  Assume it holds for $n-1$, where $n \ge 2$.  Then $x_n^4$ is the leading coefficient of $x_{n-1}$ in $f(x_{n-1},x_n)$, so by (11) and the definition of the resultant, we have that
$$R^{(n)}(x,x_n)= x_n^{4^n} \prod_{i=1}^4{R^{(n-1)}(x,\beta_i)}=\prod_{i=1}^4{x_n^{4^{n-1}}R^{(n-1)}(x,\beta_i)},$$
where $x_{n-1}=\beta_i, 1 \le i \le 4$, are the roots of the equation $f(x_{n-1},x_n)=0$.  Dividing this equation by $x_{n}^4$ and expanding with $x_{n-1}=\beta_i$ shows that
$$\beta_i^4-\left(4-\frac{8}{x_n^4}\right)\beta_i^3+6\beta_i^2-\left(4-\frac{8}{x_n^4}\right)\beta_i+1=0.$$
It follows that the elementary symmetric functions in the $\beta_i$ have degree $0$ in $x_n$, and in the product
$$R^{(n)}(x,x_n)= \prod_{i=1}^4{(x_n^{4^{n-1}}A_{n-1}(x)\beta_i^{4^{n-1}}+x_n^{4^{n-1}}S_{n-1}(x,\beta_i))},$$
the leading term is $x_n^{4^n}A_{n-1}(x)^4(\beta_1 \beta_2 \beta_3 \beta_4)^{4^{n-1}}=x_n^{4^n}A_{n-1}(x)^4$, since the product of the $\beta_i$ is $1$.  By the inductive hypothesis, the degree in $x$ of $S_{n-1}(x,x_{n-1})$ is at most $4^{n-1}-1$, so in multiplying out the remaining terms have degree at most $3 \cdot 4^{n-1}+4^{n-1}-1=4^n-1$ in $x$.  In collecting the remaining terms that involve $x_n^{4^n}$, and adding them to $A_{n-1}(x)^4$, the highest degree term in $x$ occurs only in the leading term and $A_n(x)$ is therefore monic of degree $4^n$.  It is also clear that in the product, the degrees of the terms involving $x_n$ will all be multiples of 4.  This proves part a) of the lemma.  Part b) follows immediately from a) on setting $x_n=x$.
$\square$  \bigskip

We will now show that the polynomials $R_n(x)$ have distinct roots.  \medskip

We define similar quantities for the curve
$$f_1(x,y)=\frac{f(2x,2y)}{16}=(16x^4-32x^3+24x^2-8x+1)y^4+4x^3+x.$$
We have
$$f_1(x,y) \equiv y^4+x \ (\textrm{mod} \ 2).$$
Define the iterated resultants for $f_1(x,y)$ by $\tilde R^{(1)}(x,x_1)=f_1(x,x_1)$,
\begin{align*}
\tilde R^{(2)}(x,x_2)&=Res_{x_1}(f_1(x,x_1),f_1(x_1,x_2)),\\
\tilde R^{(k)}(x,x_k)&=Res_{x_{k-1}}(\tilde R^{(k-1)}(x,x_{k-1}),f_1(x_{k-1},x_k)), \quad k \ge 3.
\end{align*}
It follows easily by induction that
$$\tilde R^{(n)}(x,x_n) \equiv x_n^{4^n}+x \ (\textrm{mod} \ 2), \ \ n\ge 1,$$
and therefore
\begin{equation}
\tilde R_n(x)=\tilde R^{(n)}(x,x) \equiv x^{4^n}+x \ (\textrm{mod} \ 2), \ \ n\ge 1.
\end{equation}
This congruence and Hensel's Lemma [h, p. 169] imply that $\tilde R_n(x)$ has at least $4^n$ distinct roots in $\textsf{K}_2$, of which $4^n-1$ are units, corresponding to the $4^n-1$ nonzero roots of the congruence (13).  Furthermore, the relation
\begin{equation}
R_n(2x)=2^{4^n} \tilde R_n(x)
\end{equation}
implies that $R_n(x)$ also has at least $4^n$ distinct roots, as well, and $N_2(k)$ monic irreducible factors of degree $k$ in $\mathbb{Z}_2[x]$, for each divisor $k$ of $2n$, where $N_2(k)$ is the number of monic irreducible polynomials of degree $k$ in $\mathbb{F}_2[x]$.  The roots $a$ of these irreducible factors (except for $a=0$, note $f(0,0)=0$) are prime elements in the ring of integers $\textsf{R}_2$ of $\textsf{K}_2$, i.e. $a \cong 2$ ($\cong$ is Hasse's notation [h3], denoting equality up to a unit factor).   \medskip

Now we make use of the identity
\begin{equation}
(x-1)^4(y-1)^4 f\left(\frac{x+1}{x-1},\frac{y+1}{y-1} \right) = 16 f(y,x).
\end{equation}
Putting
$$b=\frac{a+1}{a-1}, \ \ b_k=\frac{a_k+1}{a_k-1}, \ \ 1 \le k \le n-1,$$
where $a$ and the $a_k$ satisfy (12), the identity (15) gives that
$$f(b,b_{n-1})=0, \quad f(b_{n-1},b_{n-2})=0, \quad \cdots \quad f(b_1,b)=0.$$
It follows that $b=\frac{a+1}{a-1}$ is a root of $R_n(x)=0$ whenever $a$ is.  If $a$ is a prime element, then $b$ is clearly a unit in $\textsf{R}_2$.  This proves that $R_n(x)$ has $2\cdot 4^n$ distinct roots in $\textsf{K}_2$, for any $n \ge 1$ (including the roots $x=0, -1$), exactly half of which are units. \medskip

It follows as in [m3] that there are polynomials $\textsf{P}_n(x)$ and $\widetilde{\textsf{P}}_n(x)$ in $\mathbb{Z}[x]$ for which
\begin{equation}
R_n(x) = \prod_{k \vert n}{\textsf{P}_k(x)}, \ \ \textsf{P}_n(x)= \prod_{k|n}{R_k(x)^{\mu(n/k)}},
\end{equation}
\begin{equation}
\tilde R_n(x) = \prod_{k \vert n}{\widetilde{\textsf{P}}_k(x)},\ \ \widetilde{\textsf{P}}_n(x)= \prod_{k|n}{\tilde R_k(x)^{\mu(n/k)}},
\end{equation}
and
\begin{equation}
\textrm{deg} \ \textsf{P}_n(x) = \textrm{deg} \ \widetilde{\textsf{P}}_n(x) = 2\sum_{k \vert n}{\mu(n/k)4^k}.
\end{equation}
We note also that
\begin{equation}
R_1(x)=\textsf{P}_1(x)=x(x+1)(x^2-x+2)(x^4-4x^3+5x^2-2x+4)
\end{equation}
\begin{equation*}
\tilde R_1(x) = \widetilde{\textsf{P}}_1(x)=x(2x+1)(2x^2-x+1)(4x^4-8x^3+5x^2-x+1).
\end{equation*}
Setting
\begin{equation*}
\tilde T(z)=\frac{1}{2} T(2z), \ \ |z|_2 \le 1,
\end{equation*}
we see from (10) that
\begin{equation}
\tilde T(x) \equiv x^4 \ (\textrm{mod} \ 2), \ \ |x|_2 = 1.
\end{equation}
From (13) and (17) and the above arguments it is clear that all the irreducible factors of $\widetilde{\textsf{P}}_n(x)$ (i.e., its reduction modulo $2$) over $\mathbb{F}_4$ have degree $n$.  It is clear that $\tilde T(a)$ is a root of $\widetilde{\textsf{P}}_n(x)$ whenever the unit $a$ is, since $a$ and therefore $\tilde T(a)$ are both periodic points of $\tilde T$ with minimal period $n$.  This is because $\tilde T^k(a)=a$ for $k < n$ would imply that $a^{4^k} \equiv a$ (mod $2$), and $a$ would therefore be a root of a polynomial of degree less than $n$ over $\mathbb{F}_4$. \medskip

For such a unit $a$, $\tilde T(a)$ reduces (mod $2$) to a root of the right side of (13).  Since (13) does not have multiple roots, and by (14), half of the roots of $\widetilde{\textsf{P}}_n(x)$ are non-units, (20) shows that $a$ and $\tilde T(a)$ are roots of the same irreducible factor over $\mathbb{F}_4$, and therefore they must be roots of the same irreducible factor over $\mathbb{Q}_2$.  It follows that
$$\textsf{P}_n(x)= \prod_{i}{g_i(x) \tilde g_i(x)},$$
where the irreducible factor $g_i(x) \in \mathbb{Z}_2[x]$ has degree $n$ or $2n$;
$$\tilde g_i(x)=(x-1)^{\textrm{deg}(g_i)}g_i\left(\frac{x+1}{x-1}\right);$$
and $T$ maps the set of roots of $g_i(x)$ into itself, for each $i$.  Since $\textsf{P}_n(x) \in \mathbb{Z}[x]$, Theorem 5 implies that the minimal polynomial $b_d(x)$ of $\xi_d$ over $\mathbb{Q}$ divides $\textsf{P}_n(x)$, for any $d$ for which the automorphism $\tau_d^{-2}=\tau^{-2}$ has order $n$ in $\textrm{Gal}(\Omega_f/K)$.  In Section 5 we will prove that these are the only irreducible factors of $\textsf{P}_n(x)$, for $n>1$.

\section{A cyclic isogeny of degree $4$.}

We will now use several results from [m1, pp. 253-254] and [lm].  First, the quantity
\begin{equation*}
j_1(\alpha) = \frac{(\alpha^8-16\alpha^4+16)^3}{\alpha^8-16\alpha^4},\\
\end{equation*}
is the $j$-invariant of the elliptic curve
\begin{equation}
E_1(\alpha): \hspace{.1 in} Y^2+XY+\frac{1}{\alpha^4} Y = X^3 + \frac{1}{\alpha^4} X^2,
\end{equation}
which is the Tate normal form for a curve with a point of order $n=4$; meaning that the point $(0,0)$ has order $4$ on this curve.  Further, 
\begin{equation}
j_2(\alpha)=\frac{(\alpha^8-16\alpha^4+256)^3}{\alpha^8(\alpha^4-16)^2}
\end{equation}
is the $j$-invariant of the elliptic curve
\begin{equation*}
E_2(\alpha): \hspace{.1 in} Y^2 + XY + \frac{2}{\alpha^4} Y = X^3 + \frac{4}{\alpha^4} X^2 - \frac{1}{\alpha^8},
\end{equation*}
and $E_1(\alpha)$ is 2-isogenous to $E_2(\alpha)$ by the map $\psi_\alpha=(\psi_{\alpha,1},\psi_{\alpha,2}):E_1(\alpha) \rightarrow E_2(\alpha)$ with
$$\psi_{\alpha,1}(X)=\frac{X^2}{X+b}, \quad \psi_{\alpha,2}(X,Y)=\frac{-b^2}{X+b}+\frac{X(X+2b)Y}{(X+b)^2}, \quad b= \frac{1}{\alpha^4}.$$ \smallskip
From [lm, Eq. (4.8)] we know that $E_1(\alpha)[2]$ -- the group of $2$-torsion points on $E_1(\alpha)$ -- consists of the base point $O$, together with the points
\begin{equation}
\left(\frac{-1}{\alpha^4},0\right), \ \left(-\frac{\beta^2-4}{8\beta^2},\frac{(\beta^2-4)^2}{32\beta^4}\right), \ \left(-\frac{\beta^2+4}{8\beta^2},\frac{(\beta^2+4)^2}{32\beta^4}\right),
\end{equation}
where $16\alpha^4+16\beta^4=\alpha^4 \beta^4$.  Reversing the roles of $\alpha$ and $\beta$ in (23) gives the points of order $2$ on the curve $E_1(\beta)$. \medskip

Furthermore, still with $b = 1/\alpha^4$, the isogeny $\rho_\alpha=(\rho_{\alpha,1}(X),\rho_{\alpha,2}(X,Y))$, with
\begin{align}
\rho_{\alpha,1}(X) &= \frac{X^2-b}{X+4b}, \\
\rho_{\alpha,2}(X,Y) &= \frac{bX^2+(b-8b^2)X+3b^2-32b^3}{(X+4b)^2}+\frac{X^2+8bX+b}{(X+4b)^2} Y,
\end{align}
maps $E_2(\alpha)$ to the curve
\begin{equation}
E_3(\alpha): \hspace{.1 in} Y^2 + XY + \frac{4}{\alpha^4}Y = X^3 + \frac{16}{\alpha^4}X^2 +\frac{6}{\alpha^4}X + \frac{\alpha^4-4}{\alpha^8},
\end{equation}
and the $j$-invariant of this curve is
\begin{equation}
j_3(\alpha)=\frac{(\alpha^8-256\alpha^4+4096)^3}{\alpha^{16}(16-\alpha^4)}.
\end{equation}

We first use these facts to prove the following result.  Although we do not make explicit use of this result, we will use several of the facts mentioned in the proof in Section 5.  Moreover, the result itself is of independent interest, since it gives an interesting application for solutions of the Fermat quartic, and corresponds to the analogous result for the Fermat cubic given in [m2, Prop. 3.5].  \bigskip

\noindent {\bf Theorem 6.} {\it If $(\alpha, \beta)$ is a point on the curve
$$Fer_4: \ 16X^4+16Y^4=X^4 Y^4,$$
then there is a cyclic isogeny $\phi_{\alpha,\beta}: E_1(\alpha) \rightarrow E_1(\beta)$ of degree $4$, whose kernel is $\textrm{ker}(\phi_{\alpha,\beta})=\langle (0,0) \rangle$.}  \medskip

\noindent {\it Proof.}  The relation
$$\alpha^4=\frac{16\beta^4}{\beta^4-16}$$
implies easily using (22) that $j_2(\alpha)=j_2(\beta)$ and therefore $E_2(\alpha) \cong E_2(\beta)$.  On the other hand, there is the dual isogeny $\hat \psi_\beta: E_2(\beta) \rightarrow E_1(\beta)$.  Therefore, if $\iota: E_2(\alpha) \rightarrow E_2(\beta)$ is an isomorphism, the map
$$\phi=\hat \psi_\beta \circ \iota \circ \psi_\alpha: E_1(\alpha) \rightarrow E_1(\beta)$$
is an isogeny of degree $4$.  To determine $\textrm{ker}(\phi)$, we find an explicit isomorphism $\iota$.  Note that with $Y_1=Y+\frac{X}{2}+\frac{1}{\alpha^4}$ the equation for $E_2(\alpha)$ becomes
$$Y_1^2=X\left(X+\frac{1}{4}\right) \left(X+\frac{4}{\alpha^4}\right).$$
Using the relation
$$\frac{4}{\alpha^4}=\frac{1}{4}-\frac{4}{\beta^4}$$
and putting $X=-X_2-\frac{1}{4}, Y_1=-\sqrt{-1} Y_2$ gives the curve
\begin{equation}
Y_2^2=X_2\left(X_2+\frac{1}{4}\right) \left(X_2+\frac{4}{\beta^4}\right).
\end{equation}
Therefore, the map $\iota(X,Y)=(\iota_1(X), \iota_2(X,Y))$ can be taken to be the map
\begin{equation}
(\iota_1(X),\iota_2(X,Y))=\left(-X-\frac{1}{4}, \sqrt{-1} Y+\frac{1+\sqrt{-1}}{2} X+\frac{1+\sqrt{-1}}{\alpha^4}+\frac{1}{16} \right).
\end{equation} \smallskip

On the other hand, the $X$-coordinate of the dual isogeny $\hat \psi_\beta:E_2(\beta) \rightarrow E_1(\beta)$ is given by
\begin{equation*}
\hat \psi_{\beta,1}(X)=\frac{X^2-\frac{1}{\beta^4}}{4X+1}.
\end{equation*}
Thus, we have
$$\phi((0,0))=\hat \psi_\beta \circ \iota((0,-\frac{1}{\alpha^4}))=\hat \psi_\beta((-\frac{1}{4},\frac{1}{\alpha^4}+\frac{1}{16}))=O_1,$$
where $O_1$ is the base point on $E_1(\beta)$.  Since $\phi$ has degree $4$ and the point $(0,0)$ has order $4$, this shows that $\textrm{ker}(\phi)= \langle (0,0) \rangle $ is cyclic.  $\square$  \bigskip

We note that the $X$-coordinate of the map $\phi=\phi_{\alpha,\beta}$ is given by the rational function
\begin{equation*}
\phi_1(X)=\hat \psi_{\beta,1} \circ \iota_1 \circ \psi_{\alpha,1}(X)=\hat \psi_{\beta,1} \left(-\frac{X^2}{X+\frac{1}{\alpha^4}}-\frac{1}{4}\right)
\end{equation*}
 \begin{equation*}
=-\frac{(4\alpha^4\beta^2X^2+\alpha^4(\beta^2-4)X+\beta^2-4)(4\alpha^4\beta^2X^2+\alpha^4(\beta^2+4)X+\beta^2+4)}{64\alpha^4\beta^4X^2(\alpha^4X+1)}.
\end{equation*}
\bigskip

\eject

\section{Periodic points of $T(z)$.}

In this section we will prove the following theorem. \medskip

\noindent {\bf Theorem 7.} {\it For $n>1$, the polynomial $\textsf{P}_n(x)$ is the product of the polynomials $b_d(x)$, where $-d$ runs through all quadratic discriminants $-d \equiv 1$ (mod $8$) for which $\tau^2$ has order $n$ in the Galois group of the corresponding ring class field $\Omega_f$.  Here $\tau=\left(\frac{\Omega_f/K}{\wp_2}\right)$ is the Artin symbol (Frobenius automorphism) for the prime divisor $\wp_2$ of $2$ in $K=\mathbb{Q}(\sqrt{-d})$.}  \bigskip

Let $\xi$ be an arbitrary periodic point of $T(z)$ of minimal period $n \ge 1$ in the domain $\textsf{D}_2=\{z: 0 < |z|_2 \le \frac{1}{2}\} \subset \textsf{K}_2$, and set
\begin{equation}
\beta=2 \xi, \ \ \alpha^4=\frac{16\beta^4}{\beta^4-16}=16\frac{\xi^4}{\xi^4-1}, \ \ \beta \in \textsf{K}_2, \ \alpha \in \textsf{K}_2(\zeta_8),
\end{equation}
where $\zeta_8=\sqrt[4]{-1}$ is an eighth root of unity.  Then $(\alpha, \beta)$ is a point on $Fer_4$ (see Theorem 6) defined over $\textsf{K}_2(\zeta_8)$.  Since $\mathbb{Q}_2(\xi)$ is an unramified extension of $\mathbb{Q}_2$, and $\mathbb{Q}_2(\zeta_8)$ is totally ramified over $\mathbb{Q}_2$, there is an automorphism
\begin{equation}
\bar \tau \in \textrm{Gal}(\mathbb{Q}_2(\xi, \zeta_8)/\mathbb{Q}_2), \ \textrm{with} \  \bar \tau:=(\xi \rightarrow T(\xi), \ \zeta_8 \rightarrow \zeta_8).
\end{equation}
(Recall that $\xi$ and $T(\xi)$ are roots of the same irreducible polynomial over $\mathbb{Q}_2$, by the last assertion of Section 3.)  \smallskip

I claim now that $E_3(\beta) \cong E_1(\alpha^{\bar \tau})$, where $E_3$ and $E_1$ are the curves defined in (26) and (21).  To prove this, let $\sigma(z)$ be the linear fractional map
\begin{equation*}
\sigma(z)=\frac{2(z+2)}{z-2}.
\end{equation*}
From the fact that $f(T(\xi),\xi)=0$ we have that
\begin{equation*}
\xi^4=1-\left( \frac{T(\xi)+1}{T(\xi)-1}\right)^4
\end{equation*}
and therefore
\begin{equation*}
\left( \frac{T(\xi)+1}{T(\xi)-1}\right)^4=1-\xi^4.
\end{equation*}
Since $\beta^{\bar \tau}=2\xi^{\bar \tau}=2T(\xi)$, this gives
\begin{equation*}
\left( \frac{\beta^{\bar \tau}+2}{\beta^{\bar \tau}-2}\right)^4=1-\frac{\beta^4}{16},
\end{equation*}
and hence
\begin{equation}
\sigma(\beta^{\bar \tau})^4=16-\beta^4.
\end{equation}
Therefore, as in the proof of [lm, Prop. 8.5], and using the relation between $\alpha$ and $\beta$, we have
\begin{align*}
j(E_1(\alpha^{\bar \tau})) &= \left(\frac{(\alpha^8-16\alpha^4+16)^3}{\alpha^{4}(\alpha^4-16)}\right)^{\bar \tau}\\
&= \left(\frac{(\beta^8+224\beta^4+256)^3}{\beta^{4}(\beta^4-16)^4}\right)^{\bar \tau}\\
&=  \left(\frac{(\sigma(\beta)^8+224\sigma(\beta)^4+256)^3}{\sigma(\beta)^{4}(\sigma(\beta)^4-16)^4}\right)^{\bar \tau},
\end{align*}
since $r(z)=\frac{(z^8+224z^4+256)^3}{z^{4}(z^4-16)^4}$ is invariant under the substitution $(z \rightarrow \sigma(z))$. (See [m1, Thm. 5.2] or [lm, Section 8].)  Thus, (32) gives that
\begin{align*}
j(E_1(\alpha^{\bar \tau})) &= \frac{((16-\beta^4)^2+224(16-\beta^4)+256)^3}{(16-\beta^{4})\beta^{16}}\\
&=  \frac{(\beta^8-256\beta^4+4096)^3}{\beta^{16}(16-\beta^4)}\\
&= j(E_3(\beta)).
\end{align*}

From the isomorphism just established and the beginning remarks in Section 4, we have an isogeny
\begin{equation}
\varphi_1 =\bar \iota \circ \psi_{\alpha^{\bar \tau}} \circ \iota_3 \circ \rho_\beta
\end{equation}
of degree $4$ from $E_2(\beta)$ to $E_2(\beta^{\bar \tau})$, where $\bar \iota$ and $\iota_3$ are isomorphisms
\begin{equation*}
\bar \iota: E_2(\alpha^{\bar \tau}) \rightarrow E_2(\beta^{\bar \tau}), \ \ \iota_3: E_3(\beta) \rightarrow E_1(\alpha^{\bar \tau}).
\end{equation*}
(Note that $E_2(\alpha^{\bar \tau}) \cong E_2(\beta^{\bar \tau})$ by the beginning of the proof of Theorem 6.)  Applying the isomorphism $\bar \tau^{i-1}$ to the coefficients gives an isogeny $\varphi_i: E_2(\beta^{\bar \tau^{(i-1)}}) \rightarrow E_2(\beta^{\bar \tau^{i}})$, and therefore an isogeny
\begin{equation}
\varsigma=\varphi_n \circ \varphi_{n-1} \circ \cdots \circ \varphi_1: E_2(\beta) \rightarrow E_2(\beta),
\end{equation}
since $\bar \tau^n =1$.  This isogeny has degree $\textrm{deg}(\varsigma)=4^{n}$, and I claim that 
\begin{equation}
\Phi_{4^{n}}(j_2(\beta),j_2(\beta))=0,
\end{equation}
where $\Phi_m(X,Y)=0$ is the modular equation.  (See [co1] and [d1].)   It is well-known that (35) is equivalent to the assertion that $\textrm{ker}(\varsigma) \subset E_2(\beta)$ is cyclic.  \smallskip

From (28), the points of order $2$ on $E_2(\beta)$ are
\begin{equation}
\left(0, -\frac{1}{\beta^4}\right), \ \ \left(-\frac{1}{4}, \frac{1}{8}-\frac{1}{\beta^4}\right), \ \ \left(-\frac{4}{\beta^4}, \frac{1}{\beta^4} \right).
\end{equation}
The last of these points is in $\textrm{ker}(\rho_\beta)$, and $\rho_\beta$ maps the first two points to the point $P_1=\left( -\frac{1}{4},\frac{2}{\alpha^4} \right)$ on $E_3(\beta)$.  The other two points of order $2$ on $E_3(\beta)$ are the points
\begin{equation*}
P_2, P_3 = \left(-8\frac{\alpha^2 \pm \sqrt{-1}\beta^2}{\alpha^2 \beta^4},2 \frac{\alpha^2 \pm 2\sqrt{-1}\beta^2}{\alpha^2\beta^4} \right).
\end{equation*}
From (23), with $\alpha$ replaced by $\alpha^{\bar \tau}$, the points of order $2$ on $E_1(\alpha^{\bar \tau})$ are
\begin{equation*}
Q_1=\left(\frac{-1}{(\alpha^{\bar \tau})^4},0\right), \ Q_2=\left(-\frac{(\beta^{\bar \tau})^2-4}{8(\beta^{\bar \tau})^2},\frac{((\beta^{\bar \tau})^2-4)^2}{32(\beta^{\bar \tau})^4}\right), Q_3=\left(-\frac{(\beta^{\bar \tau})^2+4}{8(\beta^{\bar \tau})^2}, \frac{((\beta^{\bar \tau})^2+4)^2}{32(\beta^{\bar \tau})^4}\right).
\end{equation*}
\medskip

Now from (32) we have that
\begin{equation*}
\sigma(\beta^{\bar \tau})^4=-\frac{16\beta^4}{\alpha^4},
\end{equation*}
which implies that
\begin{equation*}
\sigma(\beta^{\bar \tau}) = \frac{2\beta}{\zeta_8 \alpha} ,
\end{equation*}
for some primitive eighth root of unity $\zeta_8$.  Therefore, since $\sigma$ is an involution,
\begin{equation}
\beta^{\bar \tau}= \sigma \left(\frac{2\beta}{\zeta_8\alpha}  \right) = 2 \frac{\beta+\zeta_8 \alpha}{\beta-\zeta_8\alpha}.
\end{equation}
With (37), the points of order $2$ on $E_1(\alpha^{\bar \tau})$ can be expressed in terms of $\alpha$ and $\beta$:
\begin{equation*}
Q_1=\left(-\frac{\zeta_8\alpha\beta(\beta^2+\zeta_8^2\alpha^2)}{2(\beta+\zeta_8 \alpha)^4},0\right), \ \ Q_2=\left(-\frac{\zeta_8\alpha\beta}{2(\beta+\zeta_8\alpha)^2},\frac{\zeta_8^2\alpha^2\beta^2}{2(\beta+\zeta_8\alpha)^4}\right),
\end{equation*}

\begin{equation*}
Q_3=\left(-\frac{\beta^2+\zeta_8^2\alpha^2}{4(\beta+\zeta_8\alpha)^2},\frac{(\beta^2+\zeta_8^2\alpha^2)^2}{8(\beta+\zeta_8\alpha)^4}\right).
\end{equation*}
\medskip

Converting the curves $E_3(\beta)$ and $E_1(\alpha^{\bar \tau})$ to Weierstrass normal form and using standard arguments, it can be shown that the $X$-coordinate of an isomorphism $\iota_3: E_3(\beta) \rightarrow E_1(\alpha^{\bar \tau})$ is given by
\begin{equation*}
\iota_{3,1}(X)=\frac{\beta^4+\alpha^4}{(\beta+\zeta_8\alpha)^4}X-\frac{\zeta_8\alpha(\beta^2+\zeta_8^2\alpha^2)}{2(\beta+\zeta_8\alpha)^3}.
\end{equation*}
Hence, we have that
\begin{equation*}
\iota_{3,1}\left(-\frac{1}{4}\right)=-\frac{\beta^2+\zeta_8^2\alpha^2}{4(\beta+\zeta_8\alpha)^2}.
\end{equation*}
Using (24), and comparing $X$-coordinates of the different representations of the points of order $2$ on $E_1(\alpha^{\bar \tau})$, we have
\begin{equation}
\iota_3 \circ \rho_\beta \left(0,-\frac{1}{\beta^4}\right) = Q_3 =\left(-\frac{(\beta^{\bar \tau})^2+4}{8(\beta^{\bar \tau})^2},\frac{((\beta^{\bar \tau})^2+4)^2}{32(\beta^{\bar \tau})^4}\right).
\end{equation}
Now a straightforward calculation shows that
\begin{equation}
\bar \iota_1 \circ \psi_{\alpha^{\bar \tau},1}\left(-\frac{(\beta^{\bar \tau})^2+4}{8(\beta^{\bar \tau})^2}\right)=\bar \iota_1\left(-\frac{1}{4}\right)=0,
\end{equation}
by (29), with $\alpha$ replaced by $\alpha^{\bar \tau}$.  It follows from (33), (38), (39), and (36), that
\begin{equation*}
P=\left(0,-\frac{1}{\beta^4}\right) \implies \varphi_1(P)=\left(0,-\frac{1}{(\beta^{\bar \tau})^4}\right)=P^{\bar \tau}.
\end{equation*}
Applying $\bar \tau^{i-1}$ gives that $\varphi_i(P^{\bar \tau^{i-1}})=P^{\bar \tau^{i}}$, and therefore (34) gives that
\begin{equation*}
\varsigma(P)=P^{\bar \tau^n}=P.
\end{equation*}
Since $P$ has order $2$ on $E_2(\beta)$, this shows that $P \notin \textrm{ker}(\varsigma)$.  It follows that $\textrm{ker}(\varsigma)$ is a cyclic group, and this implies (35). \smallskip

Now by a classical result [co1, p.287] we have the factorization
$$\Phi_{4^n}(x,x)=c_n \prod_{-d}{H_{-d}(x)^{r(d,4^n)}},$$
where the product is over discriminants of orders $\textsf{R}_{-d}$ of imaginary quadratic fields and
$$r(d,m)=|\{\lambda \in \textsf{R}_{-d}: \lambda \ \textrm{primitive}, \ N(\lambda)=m \}/\textsf{R}_{-d}^\times|.$$ 
The exponent $r(d,4^n)$ can only be nonzero when $4^k\cdot 4^n=x^2+dy^2$ has a primitive solution ($k=0$ or $1$).  Since $\mathbb{Q}_2(\beta)=\mathbb{Q}_2(\xi)$ is unramified and normal over $\mathbb{Q}_2$, equation (35) implies $j_2(\beta)=j(E_2(\beta))$ is a root of $H_{-d}(x)$ for some odd integer $d$; hence, $(2,xyd)=1$ and for $n>1$ we have $-d \equiv 1$ (mod 8).  \smallskip

Consequently, equation (22) shows that $\xi^4=\beta^4/16$ is a root of the polynomial
\begin{equation*}
L_d(x)=(x^2-x)^{2h(-d)}H_{-d}\left(\frac{2^8(x^2-x+1)^3}{x^2(x-1)^2}\right).
\end{equation*}
By the proof of [lm, Prop. 8.4], this polynomial factors into a product of three irreducible polynomials of degree $2h(-d)$, exactly one of which has roots which are integral for the prime $2$.  If this factor is $g(x)$, then from [lm, Eq. (8.4)] and $\textrm{deg}(g(x))=2h(-d)$ it follows that
\begin{equation}
g(x^4)=b_d(x)b_d(-x)h(x),
\end{equation}
where the irreducible polynomial $h(x)=b_d(ix)b_d(-ix)$ belongs to an extension of $\mathbb{Q}$ which is ramified over $p=2$.  Thus, $\xi$ is a root of one of the first two factors in (40).  Now the set of roots of $b_d(x)$ is stabilized by the map $\left(x \rightarrow \frac{x+1}{x-1}\right)$, and that of $b_d(-x)$ is stabilized by $\left(x \rightarrow \frac{1-x}{1+x}\right)$ (see [lm, Prop. 8.2]).  But by the factorization of $\textsf{P}_n(x)$ in Section 3, the roots of $\textsf{P}_n(x)$ are stabilized by $\left(x \rightarrow \frac{x+1}{x-1}\right)$.  If $\frac{1-\xi}{1+\xi}$ were a root of $\textsf{P}_n(x)$, then
\begin{equation*}
\frac{\frac{1-\xi}{1+\xi}+1}{\frac{1-\xi}{1+\xi}-1}=\frac{-1}{\xi}
\end{equation*}
would also be a root of $\textsf{P}_n(x)$.  But $\xi \in \textsf{D}_2$, so $-1/\xi$ is not an algebraic integer, and therefore cannot be a root of $\textsf{P}_n(x)$.  This proves that $\xi$ is a root of the polynomial $b_d(x)$ and hence that $b_d(x)$ divides $\textsf{P}_n(x)$.  From Theorem 4 and (31) we have finally that $\bar \tau=\tau^{-2}$, and since $\xi$ generates the ring class field $\Omega_f$ over $\mathbb{Q}$ and $\tau^{-2n}(\xi)=T^n(\xi)=\xi$, the automorphism $\tau^{-2}$ has order $n$ in $\textrm{Gal}(\Omega_f/K)$, where $K=\mathbb{Q}(\sqrt{-d})$.  Recalling the final remark of Section 3, this completes the proof of Theorem 7. $\square$   \bigskip

For $n=1$, we have the factorization $\textsf{P}_1(x)=x(x+1)b_7(x)b_{15}(x)$, by (19).  Hence, Theorem 7 and the formulas in (16) imply part (b) of Conjecture 1: all but two of the periodic points of $T$ in $\overline{\mathbb{Q}}_p$ generate ring class fields over $\mathbb{Q}$.  In addition, this proves Theorem 2 of the introduction, since the formulas in (16) hold over $\mathbb{Q}$, and therefore also over $\mathbb{C}$.\medskip

Denote the set of discrimimants $-d$ referred in Theorem 7 by $\mathfrak{D}_n$. Using (18) and the fact that $\textrm{deg}(b_d(x))=2h(-d)$, Theorem 7 implies the following class number relation.  \bigskip

\noindent {\bf Theorem 8.}  {\it If $h(-d)$ is the class number of the order $\textsf{R}_{-d}$ of discriminant $-d \equiv 1$ (mod $8$) in $K=\mathbb{Q}(\sqrt{-d})$, then
\begin{equation*}
\sum_{-d \in \mathfrak{D}_n}{h(-d)}=nN_4(n)=\sum_{k \vert n}{\mu (n/k)2^{2k}}, \quad n > 1,
\end{equation*}
where $\mathfrak{D}_n$ is the set of discriminants $-d \equiv 1$ (mod $8$) for which $\tau^2=\left(\frac{\Omega_f/K}{\wp_2}\right)^2$ has order $n$ in the Galois group of the corresponding ring class field $\Omega_f$.  This equation gives the total number of periodic points of $T(z)$ having minimal period $n$ in the domain $\textsf{D}_2:=\{y : 0< |y|_2 \le \frac{1}{2}\} \subset \textsf{K}_2$.  All of these periodic points (for $n>1$) are prime elements in the local field $\textsf{K}_2$.}  \bigskip

Finally, Theorem 1 summarizes the results in Proposition 3 and Theorems 4, 5, 7, and 8.

\section{Examples.}

The iterated resultants considered in Section 3 are useful in computing the polynomials $b_d(x)$ which are the minimal polynomials of the periodic points of $T(z)$.  For example, factoring $R_2(x)$ on Maple yields the polynomial $\textsf{P}_1(x)$ in (19) times
\begin{align*}
\textsf{P}_2(x)&=(x^8+20x^7+110x^6-100x^5+49x^4-80x^3-40x^2+40x+16) \\
&\times (x^8+6x^7+78x^6-84x^5+53x^4-66x^3-12x^2+24x+16) \\
&\times (x^8-6x^7+42x^6-60x^5+53x^4-54x^3+24x^2+16) \\
&= b_{63}(x)b_{55}(x) b_{39}(x).
\end{align*}
(See [lm, Section 12, Table 3].)  In addition, factoring $R_3(x)$ on Maple gives $\textsf{P}_1(x)$ times the polynomial
\begin{equation*}
\textsf{P}_3(x)=A_6(x) A_{12}(x) A_{24}(x),
\end{equation*}
where
\begin{align*}
A_6(x) &=(x^6+x^5+9x^4-13x^3+18x^2-16x+8) \\
&\times (x^6+7x^5+11x^4-15x^3+16x^2-20x+8) \\
&= b_{23}(x) b_{31}(x);
\end{align*}

\begin{align*}
A_{12}(x) &= (x^{12}-262x^{11}+20035x^{10}+13096x^9-13397x^8-15878x^7-24435x^6\\&\ -14516x^5+14372x^4+15128x^3+5440x^2+416x+64)\\
&\times (x^{12}-36x^{11}+2271x^{10}+1586x^9-1689x^8-1800x^7-2527x^6\\
&\ -2310x^5+2664x^4+832x^3+1296x^2-288x+64)\\
&\times (x^{12}-166x^{11}+8027x^{10}+5200x^9-5565x^8-6446x^7-9659x^6\\
&\ -6172x^5+6540x^4+5600x^3+2672x^2-32x+64)\\
&\times (x^{12}+16x^{11}+395x^{10}+398x^9-357x^8-316x^7-155x^6\\
&\ -1058x^5+1332x^4-704x^3+800x^2-352x+64)\\
&\times (x^{12}+184x^{11}+57491x^{10}+39206x^9-36669x^8-44260x^7-70067x^6\\
&\ -41690x^5+37644x^4+43072x^3+13616x^2+1472x+64) \\
&= b_{207}(x) b_{135}(x) b_{175}(x) b_{87}(x) b_{247}(x);
\end{align*}
and
\begin{align*}
A_{24}(x) &= (x^{24}-160x^{23}+39806x^{22}-404188x^{21}+1735295x^{20}-4082916x^{19}\\
&\ +6591016x^{18}-7995792x^{17}+7025423x^{16}-3646952x^{15}-2986282x^{14}\\
&\ +8218276x^{13}-7410127x^{12}+8124428x^{11}-590812x^{10}-4737592x^9\\
&\ +2208800x^8-5462688x^7+644992x^6+672768x^5+631808x^4+875008x^3\\
&\ +496640x^2+53248x+4096)\\
&\times (x^{24}+484x^{23}+67682x^{22}-315500x^{21}+1778351x^{20}-3320880x^{19}\\
&\ +7580476x^{18}-12603888x^{17}+15479855x^{16}-14728444x^{15}+4226978x^{14}\\
&\ +12258548x^{13}-20944063x^{12}+22569256x^{11}-11161888x^{10}-5859992x^9\\
&\ +9241280x^8-9494496x^7+2773504x^6+2227200x^5-1364224x^4+780800x^3\\
&\ +708608x^2+100352x+4096)\\
&= b_{231}(x) b_{255}(x).
\end{align*}
That each of the above polynomials is given by the corresponding $b_d(x)$ can be verified by factoring the polynomial modulo primes of the form $q=x^2+dy^2$, checking that it splits completely into linear factors (mod $q$).  Thus, we have the factorization
\begin{equation*}
\textsf{P}_3(x) = b_{23}(x) b_{31}(x) b_{207}(x) b_{135}(x) b_{175}(x) b_{87}(x) b_{247}(x) b_{231}(x) b_{255}(x)
\end{equation*}
for the periodic points of minimal period $3$.

\bigskip

\section{References.}

\noindent [ch] N. Childress, {\it Class Field Theory}, Springer, 2009. \medskip

\noindent [cn] H. Cohn, Iterated Ring Class Fields and the Icosahedron, Mathematische Annalen 255 (1981), 107-122.  \medskip

\noindent [co1] David A.Cox, {\it Primes of the Form $x^2+ny^2$; Fermat, Class Field Theory, and Complex Multiplication}, John Wiley \& Sons, 1989. \medskip

\noindent [co2] David A. Cox, {\it Galois Theory}, John Wiley \& Sons, 2004. \medskip

\noindent [d1] M. Deuring, Die Typen der Multiplikatorenringe elliptischer Funktionenk\" orper, Abh. Math. Sem. Hamb. 14 (1941), 197-272.  \medskip

\noindent [d2] M. Deuring, Die Anzahl der Typen von Maximalordnungen einer definiten Quaternionenalgebra mit primer Grundzahl, Jahresber. Deutsch. Math. Verein. 54 (1944), 24-41. \medskip

\noindent [d3] M. Deuring, Die Klassenk\"orper der komplexen Multiplikation, Enzyklop\"adie der math. Wissenschaften I2, 23 (1958), 1-60. \medskip

\noindent [h1] H. Hasse, Neue Begr\"undung der komplexen Multiplikation. I. Einordnung in die allgemeine Klassenk\"orpertheorie, J. reine angew. Math. 157 (1927), 115-139; paper 33 in {\it Mathematische Abhandlungen}, Bd. 2, Walter de Gruyter, Berlin, 1975, pp. 3-27. \medskip

\noindent [h2] H. Hasse, {\it Vorlesungen \"uber Klassenk\"orpertheorie}, Physica Verlag, W\"urzburg, 1967. \medskip

\noindent [h3] H. Hasse, {\it Number Theory}, Springer Verlag, 2002.  \medskip

\noindent [lm] R. Lynch and P. Morton, The quartic Fermat equation in Hilbert class fields of imaginary quadratic fields, http://arxiv.org/abs/1410.3008, Intl. J. of Number Theory 11 (2015), 1961-2017. \medskip

\noindent [m1] P. Morton, Explicit identities for invariants of elliptic curves, J. Number Theory 120 (2006), 234-271. \medskip

\noindent [m2] P. Morton, The cubic Fermat equation and complex multiplication on the Deuring normal form, Ramanujan J. of Math. 25 (2011), 247-275.  \medskip

\noindent [m3] P. Morton, Solutions of the cubic Fermat equation in ring class fields of imaginary quadratic fields (as periodic points of a $3$-adic algebraic function), http://arxiv.org /abs/1410.6798, to appear in Intl. J. Number Theory.  \medskip

\noindent [sch] R. Schertz, {\it Complex Multiplication}, New Mathematical Monographs, vol. 15, Cambridge University Press, 2010. \medskip

\noindent [si] J.H. Silverman, {\it Advanced Topics in the Arithmetic of Elliptic Curves}, in: Graduate Texts in Mathematics, vol. 151, Springer, New York, 1994. \bigskip

\noindent Department of Mathematical Sciences \smallskip

\noindent Indiana University - Purdue University at Indianapolis \smallskip
 
\noindent 402 N. Blackford St., LD 270 \smallskip

\noindent Indianapolis, Indiana, 46202 \smallskip

\noindent {\it e-mail: pmorton@math.iupui.edu}

\bigskip

\end{document}